# RIGIDITY OF SPACELIKE HYPERSURFACES IN SPATIALLY WEIGHTED GENERALIZED ROBERTSON-WALKER SPACETIMES

ALMA L. ALBUJER[1,*], HENRIQUE F. DE LIMA[2], ARLANDSON M. OLIVEIRA[2]
AND MARCO ANTONIO L. VELÁSQUEZ[2]

ABSTRACT. Our purpose in this paper is to apply some maximum principles in order to study the rigidity of complete spacelike hypersurfaces immersed in a spatially weighted generalized Robertson-Walker (GRW) spacetime, which is supposed to obey the so called strong null convergence condition. Under natural constraints on the weight function and on the $f$-mean curvature, we establish sufficient conditions to guarantee that such a hypersurface must be a slice of the ambient space. In this setting, we also obtain new Calabi-Bernstein type results concerning entire graphs in a spatially weighted GRW spacetime.

## 1. INTRODUCTION

The study of the geometry of spacelike hypersurfaces in certain Lorentzian ambient spaces has been a topic of growing interest in the last decades. In 1970, Calabi [14] stated the well known Calabi-Bernstein theorem: the only maximal surfaces in the 3-dimensional Lorentz-Minkowski space, that is, spacelike surfaces with zero mean curvature, are the spacelike planes. The non-parametric version of this theorem asserts that the only entire vertical graphs in $\mathbb{L}^3$ are the spacelike planes. Cheng and Yau [18] extended this result to maximal hypersurfaces in $\mathbb{L}^{n+1}$. This unicity result is no longer true for the case of spacelike hypersurfaces with non necessarily zero constant mean curvature. In fact, hyperbolic spaces are examples of non trivial spacelike hypersurfaces in $\mathbb{L}^{n+1}$ with non zero constant mean curvature. However, the conclusion is still true under some extra assumptions. Xin [32] and Aiyama [1] proved simultaneously and independently that the only spacelike hypersurfaces with constant mean curvature in $\mathbb{L}^{n+1}$ and bounded hyperbolic angle are the spacelike hyperplanes. Later on, Aledo and Alías [5] characterized spacelike hyperplanes in $\mathbb{L}^{n+1}$ as the only spacelike hypersurfaces with constant mean curvature bounded between two parallel spacelike hyperplanes.

During the last years, many authors have considered similar rigidity results in more general ambient spaces. A natural generalization of the Lorentz-Minkowski space is the class of Lorentzian product manifolds of the form $-I \times \mathbb{P}^n$, where $\mathbb{P}^n$ is an $n$-dimensional Riemannian manifold and $I \subseteq \mathbb{R}$ is an open interval. Several authors have obtained Calabi-Bernstein type results in such spaces for maximal hypersurfaces, or more generally spacelike hypersurfaces with constant mean curvature (see, for instance, Albujer and Alías [2], Albujer, Camargo and de Lima [3] and Li and Salavessa [24]). These Lorentzian product spaces are a particular case of the following family of Lorentzian manifolds: let us consider the product manifold $I \times \mathbb{P}^n$ endowed with the Lorentzian metric

$$\langle , \rangle = -dt^2 + \varrho^2 \langle , \rangle_{\mathbb{P}},$$

where $\varrho : I \to \mathbb{R}$ is a smooth positive function. The resulting Lorentzian manifold is kwown in the mathematical literature as a generalized Robertson-Walker (GRW) spacetime. The first results in this context were obtained by Alías, Romero and Sánchez in [9] for the case of compact spacelike hypersurfaces with constant mean curvature. From them on many authors have obtained unicity results not only for compact, but for complete (non-necessarily compact) spacelike hypersurfaces with constant mean curvature in a GRW satisfying certain energy conditions (see, for instance, [4, 8, 12, 13, 15, 28, 29, 30]).

Going a step further, one can consider manifolds endowed with a measure that has a smooth positive density with respect to the canonical one. The resulting spaces are the so called weighted manifolds. More precisely, given a complete pseudo-Riemannian manifold $(M, \langle , \rangle)$ and a smooth function $f : M \to \mathbb{R}$, the weighted manifold $M_f$ associated to $M$ and $f$ is the triple $(M, \langle , \rangle, d\mu = e^{-f} dM)$, where $dM$ is the standard volume element of $M$. In this context, Bakry and Émery introduced in [10] the Bakry-Émery-Ricci tensor as a suitable generalization of the classical Ricci tensor, Ric, defined by

$$\mathrm{Ric}_f = \mathrm{Ric} + \mathrm{Hess}\, f.$$

It seems natural to try to extend results stated in terms of the Ricci curvature tensor to analogous results for the Bakry-Émery-Ricci tensor.

It is interesting to remark that weighted manifolds are closely related to some classical mathematical concepts, as they can be used as a powerful mathematical tool in order to obtain new results related to them. Specifically, in the case where $\mathrm{Ric}_f$ is constant, we can induce on $M$ a structure of a gradient Ricci soliton. Its mathematical relevance is due to the Perelman's solution of the Poincaré conjecture since gradient Ricci solitons correspond to self-similar solutions to the Hamilton's Ricci flow and often arise as limits of dilations of singularities developed along the Ricci flow. For an overview of results in this context one can consult [27]. On the other hand, weighted manifolds have also been considered when studying harmonic heat flows and heat kernels. For instance, Grigor'yan and Saloff-Coste established in [21] a result which relates the heat kernel on a complete, noncompact Riemannian manifold $M$ with the Dirichlet heat kernel on the exterior of a compact set of $M$. They got a relation between both heat kernels either in the case when $M$ is transient or when it is recurrent, and this was possible since they applied in the proof the theory of weighted manifolds.

In this paper, we consider complete spacelike hypersurfaces in a spatially weighted GRW spacetime obeying a suitable energy condition, and under some extra assumptions on the geometry of the hypersurface we conclude several rigidity results. The manuscript is organized as follows: in Section 2, we introduce some basic notions and facts related to spacelike hypersurfaces in GRW spacetimes and weighted manifolds in general. In Section 3, we present some auxiliary lemmas that will be needed in the proofs of our main results. In particular we state a weak Omori-Yau maximum principle for the drifted Laplacian (Lemma 3) and some other results that will be crucial in order to obtain our main results which will be presented in the last two sections. Theorems 1 and 2 in Section 4 are obtained as a consequence of the weak Omori-Yau maximum principle for the drifted Laplacian. It is interesting to observe that Theorem 1 is a refinement of [17, Theorem 2], while Theorem 2 is a natural generalization of [4, Theorem 3.3]. Our last result in this section (Theorem 3) is a completely new rigidity result obtained from Lemmas 4 and 5. In Section 5, we recall the concept of entire vertical graph in a spatially weighted GRW spacetime and we give non-parametric versions of the above results (Theorems 4, 5 and 6).

## 2. Set up

Let $(\mathbb{P}^n, \langle , \rangle_{\mathbb{P}})$ be an $n$-dimensional, connected, oriented Riemannian manifold, $I \subseteq \mathbb{R}$ an open interval and $\varrho : I \to \mathbb{R}$ a positive smooth function. Let us consider the product manifold $I \times \mathbb{P}^n$, and denote by $\pi_I$ and $\pi_{\mathbb{P}}$ the projections onto the factors $I$ and $\mathbb{P}^n$, respectively. The class of Lorentzian manifolds which will be of our concern here is the one obtained by furnishing $I \times \mathbb{P}^n$ with the Lorentzian metric

$$\langle v, w \rangle_p = -\pi_I^*(v)\pi_I^*(w) + (\varrho \circ \pi_I)(p)^2 \langle \pi_{\mathbb{P}}^*(v), \pi_{\mathbb{P}}^*(w) \rangle_{\mathbb{P}_{\pi_{\mathbb{P}}(p)}},$$

for all $p \in I \times \mathbb{P}^n$ and $v, w \in T_p(I \times \mathbb{P}^n)$. In such a case, we write $-I \times_\varrho \mathbb{P}^n$. Let us observe that $-I \times_\varrho \mathbb{P}^n$ is a Lorentzian *warped product* with *warping function* $\varrho$ and fiber $\mathbb{P}^n$. When $\mathbb{P}^n$ has constant sectional curvature, the warped product $-I \times_\varrho \mathbb{P}^n$ has been known in the mathematical literature as a Robertson-Walker (RW) spacetime, an allusion to the fact that, for $n = 3$, it is an exact solution of the Einstein's field equations (cf. Chapter 12 of [25]). After [9], such warped products have usually been referred to as *generalized Robertson-Walker (GRW) spacetimes*, and we shall stick to this usage along this paper.

Let $\Sigma^n$ be an $n$-dimensional connected manifold. A smooth immersion $\psi : \Sigma^n \to -I \times_\varrho \mathbb{P}^n$ is said to be a *spacelike hypersurface* if $\Sigma^n$, furnished with the metric induced from $\langle , \rangle$ via $\psi$, is a Riemannian manifold. If this is so, we shall always assume that the metric on $\Sigma^n$ is the induced one, which will also be denoted by $\langle , \rangle$. In this setting, it follows from the connectedness of $\Sigma^n$ that one can uniquely choose a globally defined timelike unit normal vector field $N \in \mathfrak{X}(\Sigma)^\perp$, having the same time-orientation of $\partial_t$, i.e., such that $\langle N, \partial_t \rangle \leq -1$. One then says that $N$ is the *future-pointing Gauss map* of $\Sigma^n$.

Now, we consider two particular functions naturally attached to a spacelike hypersurface $\Sigma^n$ immersed into a GRW spacetime $-I \times_\varrho \mathbb{P}^n$, namely, the (vertical) *height function* $h = (\pi_I)|_\Sigma$ and the *angle function* $\eta = \langle N, \partial_t \rangle$. Given any vector field $V \in \mathfrak{X}(-I \times_\varrho \mathbb{P}^n)$, we denote by $V^* = \pi_\mathbb{P}^*(V) = V + \langle V, \partial_t \rangle \partial_t$ the projection of $V$ onto $\mathbb{P}^n$. In particular, $N^* = N + \eta \partial_t$, and therefore

$$\langle N^*, N^* \rangle = \eta^2 - 1. \tag{2.1}$$

Let us denote by $\overline{\nabla}$ and $\nabla$ the gradients with respect to the metrics of $-I \times_\varrho \mathbb{P}^n$ and $\Sigma^n$, respectively. Then, a simple computation shows that the gradient of $\pi_I$ on $-I \times_\varrho \mathbb{P}^n$ is given by

$$\overline{\nabla} \pi_I = -\langle \overline{\nabla} \pi_I, \partial_t \rangle \partial_t = -\partial_t,$$

so that the gradient of $h$ on $\Sigma^n$ is

$$\nabla h = (\overline{\nabla} \pi_I)^\top = -\partial_t^\top = -\partial_t - \langle N, \partial_t \rangle N. \tag{2.2}$$

Thus, from (2.2) we get

$$|\nabla h|^2 = \eta^2 - 1, \tag{2.3}$$

where $|\ |$ denotes the norm of a vector field on $\Sigma^n$.

At this point we recall that, given a complete semi-Riemannian manifold $(M, \langle,\rangle)$ and a smooth function $f$ on $M$, the *weighted manifold* $M_f$ associated to $M$ and $f$ is the triple $(M, \langle,\rangle, d\mu = e^{-f} dM)$, where $dM$ is the volume element of $M$. In this setting, we will consider the so-called Bakry-Émery-Ricci tensor, introduced by Bakry and Émery in [10] as a suitable extension of the standard Ricci tensor, Ric, which is defined by

$$\mathrm{Ric}_f = \mathrm{Ric} + \mathrm{Hess}\, f. \tag{2.4}$$

From a splitting theorem due to Case (cf. [16, Theorem 1.2]) it follows that if $-I \times_\varrho \mathbb{P}^n$ is a weighted GRW spacetime endowed with a weight function $f$ which is bounded from above and such that $\overline{\mathrm{Ric}}_f(V, V) \geq 0$ for all timelike vector fields $V$, then $f$ must be constant along $\mathbb{R}$. Motivated by this result, along this work we will consider *spatially weighted GRW spacetimes* $-I \times_\varrho \mathbb{P}^n_f$, which means that the weight function $f$ does not depend on the parameter $t \in I$ or, in other words, $\langle \overline{\nabla} f, \partial_t \rangle = 0$.

In this setting, for a spacelike hypersurface $\Sigma^n$ immersed in $-I \times_\varrho \mathbb{P}^n$, the *f-divergence operator* on $\Sigma^n$ is defined by

$$\mathrm{div}_f(X) = e^f \mathrm{div}(e^{-f} X),$$

for all tangent vector fields $X$ on $\Sigma^n$ and, given a smooth function $u : \Sigma^n \to \mathbb{R}$, its *drifted Laplacian* is defined by

$$\Delta_f u = \mathrm{div}_f(\nabla u) = \Delta u - \langle \nabla u, \nabla f \rangle. \tag{2.5}$$

Finally, according to Gromov [22] the *f-mean curvature* $H_f$ of $\Sigma^n$ is defined by

$$n H_f = n H - \langle \overline{\nabla} f, N \rangle, \tag{2.6}$$

where $H$ denotes the standard mean curvature of $\Sigma^n$ with respect to its future-pointing Gauss map $N$. In analogy to the case of the standard mean curvature, the $f$-mean curvature $H_f$ on $\Sigma^n$ is related to the variational problem for the weighted area functional

$$vol_f(\Sigma) = \int_\Sigma e^{-f} d\Sigma.$$

## 3. Auxiliary results

In order to prove our rigidity results in spatially weighted GRW spacetimes of the type $-I \times_\varrho \mathbb{P}^n_f$, we will need some auxiliary lemmas. In the first two lemmas, we present some useful computations that will be essential in the proofs of our main results.

**Lemma 1.** *Let $\Sigma^n$ be a spacelike hypersurface immersed in a spatially weighted GRW spacetime $-I \times_\varrho \mathbb{P}^n_f$, with height function $h$ and angle function $\eta$. Then,*

$$\begin{aligned}\Delta_f(\varrho(h)\eta) =&\, n\varrho(h)\langle \nabla H_f, \partial_t \rangle + n\varrho'(h) H_f + \varrho(h)\eta |A|^2 + \varrho(h)\eta \overline{\mathrm{Hess}} f(N, N) \\ &+ \varrho(h)\eta \left( \widetilde{\mathrm{Ric}}(N^*, N^*) - (n-1)(\log \varrho)''(h)|\nabla h|^2 \right),\end{aligned} \tag{3.1}$$

where $A$ is the shape operator of $\Sigma^n$ related to $N$, $|A|$ denotes its norm, $\sigma(t) = \int_{t_0}^{t} \varrho(s)ds$ and $\widetilde{\mathrm{Ric}}$ stands for the Ricci curvature tensor of $\mathbb{P}^n$.

*Proof.* In [6, Corollary 8.2] it is proven that

$$\begin{aligned}
(3.2) \quad \Delta(\varrho(h)\eta) =& n\varrho(h)\langle \nabla H, \partial_t\rangle + n\varrho'(h)H + \varrho(h)\eta |A|^2 \\
& + \varrho(h)\eta \left( \widetilde{\mathrm{Ric}}(N^*, N^*) - (n-1)(\log \varrho)''(h)|\nabla h|^2 \right).
\end{aligned}$$

Taking into acccount (2.6), it follows that

$$n\varrho(h)\langle \nabla H, \partial_t \rangle = n\varrho(h)\langle \nabla H_f, \partial_t \rangle + \varrho(h) \partial_t^\top (\langle \overline{\nabla} f, N\rangle).$$

Moreover, from a straightforward computation we get

$$\partial_t^\top(\langle \overline{\nabla} f, N\rangle) = -\frac{\varrho'}{\varrho}(h)\langle \overline{\nabla} f, N\rangle + \eta \, \overline{\mathrm{Hess}}\, f(N, N) - \langle \overline{\nabla} f, A\partial_t^\top \rangle$$

and

$$\nabla(\varrho(h)\eta) = -\varrho(h) A \partial_t^\top.$$

So (3.2) can be written as

$$\begin{aligned}
(3.3) \quad \Delta(\varrho(h)\eta) =& n\varrho(h)\langle \nabla H_f, \partial_t\rangle - \varrho'(h)\langle \overline{\nabla} f, N\rangle + \varrho(h)\eta \overline{\mathrm{Hess}}\, f(N, N) \\
& + \langle \overline{\nabla} f, \nabla(\varrho(h)\eta)\rangle + n\varrho'(h)H + \varrho(h)\eta |A|^2 \\
& + \varrho(h)\eta \left( \widetilde{\mathrm{Ric}}(N^*, N^*) - (n-1)(\log \varrho)''(h)|\nabla h|^2 \right).
\end{aligned}$$

Finally, (3.1) follows from (3.3) and (2.5). $\square$

In what follows, a slab $[t_1, t_2] \times \mathbb{P}^n = \{(t,q) \in -I \times_\varrho \mathbb{P}^n : t_1 \leq t \leq t_2\}$ is called a *timelike bounded region* of the spatially weighted GRW spacetime $-I \times_\varrho \mathbb{P}_f^n$. Our second auxiliary result gives sufficient conditions to guarantee that the Bakry-Émery-Ricci curvature tensor of a spacelike hypersurface immersed in a spatially weighted GRW spacetime is bounded from below. According to the terminology established by Alías and Colares in [6], we say that a spatially weighted GRW spacetime $-I \times_\varrho \mathbb{P}_f^n$ obeys the *strong null convergence condition* (SNCC) when the sectional curvatures $K_\mathbb{P}$ of its Riemannian fiber $\mathbb{P}$ satisfy the following inequality

$$(3.4) \quad K_\mathbb{P} \geq \sup_I (\varrho^2 (\log \varrho)'').$$

The following lemma establishes some sufficient conditions in order to guarantee that the Bakry-Émery-Ricci tensor is bounded from below. In [19, Lemma 3] we can find a similar result where different hypothesis have been considered.

**Lemma 2.** *Let $-I \times_\varrho \mathbb{P}^n$ be a spatially weighted GRW spacetime obeying (3.4) and such that the Hessian of the weight function $f$ is bounded from below. Let $\Sigma^n$ be a spacelike hypersurface which lies in a timelike bounded region of $-I \times_\varrho \mathbb{P}_f^n$. Suppose that the $f$-mean curvature $H_f$ is bounded on $\Sigma^n$, then the Bakry-Émery-Ricci tensor $\mathrm{Ric}_f$ of $\Sigma^n$ is bounded from below.*

*Proof.* We recall that the curvature tensor $R$ of a spacelike hypersurface $\psi : \Sigma^n \to -I \times_\varrho \mathbb{P}_f^n$ can be described in terms of the shape operator $A$ and the curvature tensor $\overline{R}$ of $-I \times_\varrho \mathbb{P}_f^n$ by the so-called Gauss equation given by

$$(3.5) \quad R(X, Y)Z = (\overline{R}(X, Y)Z)^\top - \langle AX, Z\rangle AY + \langle AY, Z\rangle AX,$$

for every tangent vector fields $X, Y, Z \in \mathfrak{X}(\Sigma)$. Here, as in [25], the curvature tensor $R$ is given by

$$R(X, Y)Z = \nabla_{[X,Y]}Z - [\nabla_X, \nabla_Y]Z,$$

where $[\,,]$ denotes the Lie bracket and $X, Y, Z \in \mathfrak{X}(\Sigma)$.

Let us consider $X \in \mathfrak{X}(\Sigma)$ and a local orthonormal frame $\{E_1, \ldots, E_n\}$ of $\mathfrak{X}(\Sigma)$. Then, it follows from the Gauss equation (3.5) that

$$\text{Ric}(X,X) = \sum_{i=1}^{n}\langle \overline{R}(X,E_i)X,E_i\rangle + nH\langle AX,X\rangle + |AX|^2. \tag{3.6}$$

Moreover, we have that (see [25] for details)

$$\begin{aligned}\overline{R}(X,Y)Z =& \widetilde{R}(X^*,Y^*)Z^* + ((\log \varrho)'(h))^2(\langle X,Z\rangle Y - \langle Y,Z\rangle X) \\ &+ (\log \varrho)''(h)\langle Z,\partial_t\rangle(\langle Y,\partial_t\rangle X - \langle X,\partial_t\rangle Y) \\ &- (\log \varrho)''(h)(\langle Y,\partial_t\rangle\langle X,Z\rangle - \langle X,\partial_t\rangle\langle Y,Z\rangle)\partial_t,\end{aligned}$$

$\widetilde{R}$ being the curvature tensor of $\mathbb{P}^n$, and hence

$$\begin{aligned}\langle \overline{R}(X,E_i)X,E_i\rangle =& \varrho(h)^2 K_{\mathbb{P}}(X^*,E_i^*)(|X^*|_{\mathbb{P}}^2|E_i^*|_{\mathbb{P}}^2 - \langle X^*,E_i^*\rangle_{\mathbb{P}}^2) \\ &+ ((\log \varrho)'(h))^2(|X|^2 - \langle X,E_i\rangle^2) \\ &+ (\log \varrho)''(h)\langle X,\nabla h\rangle(\langle \nabla h,E_i\rangle\langle X,E_i\rangle - \langle X,\nabla h\rangle) \\ &- (\log \varrho)''(h)(\langle \nabla h,E_i\rangle|X|^2 - \langle X,\nabla h\rangle\langle X,E_i\rangle)\langle \nabla h,E_i\rangle.\end{aligned} \tag{3.7}$$

On the other hand, one can easily see that

$$\begin{aligned}|X^*|_{\mathbb{P}}^2|E_i^*|_{\mathbb{P}}^2 - \langle X^*,E_i^*\rangle_{\mathbb{P}}^2 =& \\ \frac{1}{\varrho(h)^4}(|X|^2 &+ \langle X,\nabla h\rangle^2 + |X|^2\langle \nabla h,E_i\rangle^2 - \langle X,E_i\rangle^2 \\ &- 2\langle X,\nabla h\rangle\langle X,E_i\rangle\langle \nabla h,E_i\rangle),\end{aligned}$$

which jointly with (3.4) and (3.7) implies the following lower bound

$$\sum_{i=1}^n \langle \overline{R}(X,E_i)X,E_i\rangle \geq \frac{\varrho''}{\varrho}(h)(n-1)|X|^2, \tag{3.8}$$

for all $X \in \mathfrak{X}(\Sigma)$. Thus, from (3.6) and (3.8), we get

$$\text{Ric}(X,X) \geq \frac{\varrho''}{\varrho}(h)(n-1)|X|^2 + nH\langle AX,X\rangle + |AX|^2. \tag{3.9}$$

Since the Hessian of f is bounded from below, we have

$$\begin{aligned}\text{Hess}\, f(X,X) &= \overline{\text{Hess}}\, f(X,X) - \langle \overline{\nabla} f,N\rangle\langle AX,X\rangle \\ &\geq -\beta|X|^2 - \langle \overline{\nabla} f,N\rangle\langle AX,X\rangle,\end{aligned} \tag{3.10}$$

for certain constant $\beta$. Therefore, from (2.4), (2.6), (3.9) and (3.10), we obtain

$$\text{Ric}_f(X,X) \geq \left(\frac{\varrho''}{\varrho}(h)(n-1) - \beta\right)|X|^2 + nH_f\langle AX,X\rangle + |AX|^2. \tag{3.11}$$

We can write

$$nH_f\langle AX,X\rangle + |AX|^2 = \left|AX + \frac{nH_f}{2}X\right|^2 - \frac{n^2H_f^2}{4}|X|^2,$$

so inequality (3.11) becomes

$$\text{Ric}_f(X,X) \geq \left(\frac{\varrho''}{\varrho}(h)(n-1) - \beta\right)|X|^2 + \left|AX + \frac{nH_f}{2}X\right|^2 - \frac{n^2H_f^2}{4}|X|^2 \tag{3.12}$$

for all $X \in \mathfrak{X}(\Sigma)$. Finally, the assumptions that $H_f$ is bounded and $\Sigma^n$ is contained in a timelike bounded region of $-I \times_\varrho \mathbb{P}_f^n$ guarantee that the Bakry-Émery-Ricci tensor $\text{Ric}_f$ of $\Sigma^n$ is bounded from below. □

The following key lemma is a weak Omori-Yau's generalized maximum principle for the drifted Laplacian. A proof of it can be found in [27, Remark 2.18, Chapter 2] (see also [11, Theorem 5.4]).

**Lemma 3.** *Let $(M^n, \langle,\rangle, e^{-f}dM)$ be a complete weighted manifold whose Bakry-Émery-Ricci curvature tensor is bounded from below and let $\varphi : M \to \mathbb{R}$ be a smooth function bounded from below on $M^n$. Then, there exists a sequence of points $\{p_k\}_{k\in\mathbb{N}} \subset M^n$ such that*

$$\lim_k \varphi(p_k) = \inf_M \varphi \quad \text{and} \quad \liminf_k \Delta_f \varphi(p_k) \geq 0.$$

Let us consider $\mathcal{L}_f^p(M) := \{u : M^n \to \mathbb{R} : \int_M |u|^p(x) e^{-f(x)} dM < +\infty\}$. The following result is a consequence of [26, Theorem 1.1].

**Lemma 4.** *Let $u$ be a nonnegative smooth $f$-subharmonic function on a complete Riemannian manifold $M^n$. If $u \in \mathcal{L}_f^p(M)$, for some $p > 1$, then $u$ is constant.*

Finally, we will need the next result due to Wei and Wylie [31].

**Lemma 5.** *All complete noncompact Riemannian manifolds with nonnegative Bakry-Émery-Ricci tensor for some bounded weight function $f$ have at least linear $f$-volume growth (i.e., for any $p \in \Sigma^n$ $\mathrm{vol}_f(B(p,R))$ has at least linear growth on $R$, where $B(p,R)$ is the geodesic ball in $\Sigma^n$ centered at $p$ with radius $R$).*

## 4. Uniqueness results in spatially weighted GRW spacetimes

Our first two rigidity results are a consequence of the weak Omori-Yau maximum principle for the drifted Laplacian (Lemma 3). The first of them generalizes Theorem 4.3 in [15].

**Theorem 1.** *Let $-I \times_\varrho \mathbb{P}_f^n$ be a spatially weighted GRW spacetime obeying (3.4) and such that the Hessian of the weight function $f$ is bounded from below. Let $\psi : \Sigma^n \to -I \times_\varrho \mathbb{P}_f^n$ be a complete spacelike hypersurface which lies in a timelike bounded region of $-I \times_\varrho \mathbb{P}_f^n$. Suppose that the $f$-mean curvature $H_f$ of $\Sigma^n$ satisfies*

$$(4.1) \qquad (\log \varrho)'(h) \leq H_f \leq \alpha \quad \text{and} \quad H_f \geq 0,$$

*for some positive constant $\alpha$. If*

$$(4.2) \qquad |\nabla h| \leq \beta \inf_\Sigma |H_f - (\log \varrho)'(h)|^\gamma$$

*for some constants $\beta$ and $\gamma$, $\beta > 0$, $\gamma \neq 0$, then $\Sigma^n$ is a slice $\{t\} \times \mathbb{P}^n$.*

*Proof.* From [17, Lemma 1, (ii)], for $\sigma(t) = \int_{t_0}^t \varrho(s) ds$ we get

$$\Delta_f \sigma(h) = -n\varrho(h) \left((\log \varrho)'(h) + \eta H_f\right)$$
$$\geq n\varrho(h) \left(H_f - (\log \varrho)'(h)\right).$$

Thus, by (4.1), we conclude that $\Delta_f \sigma(h) \geq 0$ on $\Sigma^n$. Since $\sigma$ is bounded from above and, by Lemma 2, $\mathrm{Ric}_f$ is bounded from below, the hypothesis of Lemma 3 are hold and we can take a sequence of points $\{p_k\}_{k\in\mathbb{N}} \subset \Sigma^n$ such that

$$0 \geq \limsup_k \Delta_f \sigma(h(p_k)) \geq \lim_k \left[n\varrho(h) \left(H_f - (\log \varrho)'(h)\right)\right](p_k) \geq 0.$$

Since $\Sigma^n$ is contained in a timelike bounded region of $-I \times_\varrho \mathbb{P}_f^n$, there exists a positive constant $C$ such that $\varrho(h(p)) \geq C$, for all $p \in \Sigma^n$. Therefore, we have that $\lim_k (H_f - (\log \varrho)'(h))(p_k) = 0$ and, taking into account our hypothesis (4.2), we conclude the proof. $\square$

**Remark 1.** *Along the preparation of this manuscript, the authors have realized that the proof of [17, Theorem 2] is not correct, since it is necessary to ask $H_f \geq 0$ in order to sign the drifted Laplacian of $\sigma(h)$. Having this into account, the above result can be seen as a refinement of [17, Theorem 2].*

**Theorem 2.** *Let $-I \times_\varrho \mathbb{P}_f^n$ be a spatially weighted GRW spacetime obeying (3.4) with convex weight function $f$ (i.e., $\overline{\mathrm{Hess}} f \geq 0$). Let $\psi : \Sigma^n \to -I \times_\varrho \mathbb{P}_f^n$ be a complete spacelike hypersurface which lies in a timelike bounded region and with constant $f$-mean curvature $H_f$ satisfying*

$$(4.3) \qquad 0 \leq H_f \sup_\Sigma (\log \varrho)'(h) \leq H^2.$$

*If*

$$(4.4) \qquad |\nabla h|^2 \leq \alpha \left( \inf_\Sigma H^2 - H_f \sup_\Sigma (\log \varrho)'(h) \right)^\beta$$

*for some constants $\alpha$ and $\beta$, $\alpha > 0$, $\beta \neq 0$ then $\Sigma^n$ is a slice $\{t\} \times \mathbb{P}^n$.*

*Proof.* From inequality (3.4) we get, in particular, that

$$\widetilde{\mathrm{Ric}}(N^*, N^*) - (n-1)(\log \varrho)''(h) |\nabla h|^2 \geq 0.$$

Therefore, from (3.1) and the assumptions of the theorem it holds

$$(4.5) \qquad \Delta_f(\varrho(h)\eta) \leq n\varrho'(h) H_f + \varrho(h)\eta |A|^2.$$

From (4.5) and (4.3) we get

$$\Delta_f(\varrho(h)\eta) \leq \varrho(h)\eta \left( |A|^2 - n H_f \sup_\Sigma (\log \varrho)'(h) \right).$$

On the other hand $|A|^2 = n^2 H^2 - n(n-1) H_2$, where $H_2$ is the second order mean curvature defined by $\binom{n}{2} H_2 = \sum_{i<j} k_i k_j$, being $k_i$, $i = 1, ..., n$ the main curvatures of $\Sigma^n$. Therefore,

$$(4.6) \begin{aligned} \Delta_f(\varrho(h)\eta) &\leq \varrho(h)\eta \left( n^2 H^2 - n(n-1) H_2 - n H_f \sup_\Sigma (\log \varrho)'(h) \right) \\ &= n(n-1)\varrho(h)\eta (H^2 - H_2) + n\varrho(h)\eta \left( H^2 - H_f \sup_\Sigma (\log \varrho)'(h) \right) \leq 0, \end{aligned}$$

where we have used again (4.3) and the fact that $H^2 - H_2 \geq 0$.

From (4.4) and (2.3) we have that $\eta$ is bounded, which jointly with the fact that $\Sigma^n$ lies in a timelike bounded region of $-I \times_\varrho \mathbb{P}_f^n$ implies that the function $\varrho(h)\eta$ is bounded from below. On the other hand, by Lemma 2 we know that $\mathrm{Ric}_f$ is bounded from below, so we can apply the weak Omori-Yau maximum principle for the drifted Laplacian (Lemma 3), and conclude that there exists a sequence of points $\{p_k\} \subset \Sigma^n$ such that

$$\lim_k \varrho(h(p_k))\eta(p_k) = \inf_\Sigma \varrho(h)\eta \qquad \text{and} \qquad \liminf_k \Delta_f \varrho(h(p_k))\eta(p_k) \geq 0.$$

Thus, (4.6) implies that

$$\begin{aligned} 0 &\leq \liminf_k \Delta_f \varrho(h(p_k))\eta(p_k) \\ &\leq n(n-1) \lim_k \left( \varrho(h(p_k))\eta(p_k)(H^2 - H_2)(p_k) \right) \\ &\quad + n \lim_k \left( \varrho(h(p_k))\eta(p_k) \left( H^2(p_k) - H_f \sup_\Sigma (\log \varrho)'(h) \right) \right) \leq 0, \end{aligned}$$

so in particular

$$\lim_k H^2(p_k) - H_f \sup_\Sigma (\log \varrho)'(h) = 0.$$

Consequently, from (4.3) we get

$$\inf_\Sigma H^2 - H_f \sup_\Sigma (\log \varrho)'(h) = 0,$$

and the result follows from (4.4). □

**Remark 2.** *Observe that in the particular case where the function $f$ is constant, the $f$-mean curvature is just the usual mean curvature $H$. In this sense, the above result is a generalization of [4, Theorem 3.3]. It is interesting to observe that although the results in [4] are referred to hypersurfaces in RW spacetimes, the same conclusions can be obtained when considering GRW spacetimes $-I \times_\varrho \mathbb{P}^n$ such that the sectional curvatures of the fiber are bounded from below, as it happens under the SNCC condition.*

Our third main result is obtained as an application of Lemmas 4 and 5.

**Theorem 3.** *Let $-I \times_\varrho \mathbb{P}_f^n$ be a spatially weighted GRW spacetime and let $\psi : \Sigma^n \to -I \times_\varrho \mathbb{P}_f^n$ be a complete spacelike hypersurface. Suppose that $H_f > 0$, $\varrho'(h) > 0$ and that the following inequalities are satisfied*

$$\frac{n^2}{4}((\log \varrho)'(h))^2 \leq \frac{n^2 H_f^2}{4} \leq (n-1)\frac{\varrho''}{\varrho}(h). \tag{4.7}$$

*If $\varrho(h) \in \mathcal{L}_f^p(\Sigma)$, for some $p > 1$, then $\Sigma^n$ is a slice of $-I \times_\varrho \mathbb{P}_f^n$ with $vol_f(\Sigma) < +\infty$. In addition, if $-I \times_\varrho \mathbb{P}_f^n$ obeys the SNCC (3.4) and $f$ is bounded and convex, then $\Sigma^n$ is compact.*

*Proof.* From [17, Lemma 1, (i)] we get

$$\Delta_f h = -(\log \varrho)'(h)(n + |\nabla h|^2) - n H_f \eta. \tag{4.8}$$

Moreover, we note that it is not difficult to verify that the hypothesis (4.7) implies that $(\log \varrho)''(h) \geq 0$ on $\Sigma^n$. Consequently, using the assumptions of the theorem, from (4.8) we get

$$\begin{aligned}\Delta_f \varrho(h) &= \varrho'(h)\Delta_f h + \varrho''(h)|\nabla h|^2 \\ &\geq n\varrho'(h)\left(H_f - (\log \varrho)'(h)\right) \geq 0.\end{aligned} \tag{4.9}$$

Thus, since we are assuming that $\varrho(h) \in \mathcal{L}_f^p(\Sigma)$, from (4.9) we can apply Lemma 4 to conclude that $\varrho(h)$ is constant on $\Sigma^n$ and $vol_f(\Sigma) < +\infty$. Hence, since $\varrho'(h) > 0$ on $\Sigma^n$, we get that $h$ is also constant and, consequently, $\Sigma^n$ must be a slice of $-I \times_\varrho \mathbb{P}_f^n$.

Furthermore, if $-I \times_\varrho \mathbb{P}_f^n$ obeys SNCC and the weight function $f$ is convex, from (3.12) we obtain that

$$\text{Ric}_f(X, X) \geq \left((n-1)\frac{\varrho''}{\varrho}(h) - \frac{n^2 H_f^2}{4}\right)|X|^2, \tag{4.10}$$

for all $X \in \mathfrak{X}(\Sigma)$. Thus, considering (4.7) into (4.10), we see that $\text{Ric}_f$ is nonnegative. Therefore, Lemma 5 guarantees that $\Sigma^n$ must be compact. □

**Remark 3.** *Let us consider the closed conformal vector field $V = V(t, p) = \varrho(t)\partial_t$ which is globally defined on the GRW spacetime $-I \times_\varrho \mathbb{P}^n$. Following the ideas of [20] (see also [23]), given a spacelike hypersurface $\Sigma^n$ immersed in $-I \times_\varrho \mathbb{P}_f^n$ with future-pointing Gauss map $N$, we can write $V_q = E(q)N_q + V_q^\top$, for each $q \in \Sigma^n$, where $E(q) := -\langle V_q, N_q \rangle > 0$ and $V_q^\top$ are, respectively, the energy and the n-momentum that the instantaneous observer $N_q$ measures for $V_q$.*

*Considering the spatially weighted GRW spacetime $-I \times_\varrho \mathbb{P}_f^n$ and extending the concept of total energy already established in [20], we say that $\Sigma^n$ has finite total $(f, p)$-energy when $\int_\Sigma E^p(q)e^{-f(q)}d\Sigma < +\infty$, for $p \geq 1$. So, since $E(q) \geq \varrho(h(q))$ for all $q \in \Sigma^n$, if we assume in Theorem 3 that $\Sigma^n$ has finite total $(f, p)$-energy instead of $\varrho(h) \in \mathcal{L}_f^p(\Sigma)$, the result still holds.*

## 5. Entire vertical graphs in a spatially weighted GRW spacetime

Let $\Omega \subseteq \mathbb{P}^n$ be a connected domain and let $u \in \mathcal{C}^\infty(\Omega)$ be a smooth function, then $\Sigma^n(u)$ will denote the vertical graph over $\Omega$ determined by $u$, that is,

$$\Sigma^n(u) = \{(u(x), x) : x \in \Omega\} \subset -I \times_\varrho \mathbb{P}^n.$$

The graph is said to be entire if $\Omega = \mathbb{P}^n$. The metric induced on $\Omega$ from the Lorentzian metric of the ambient space via $\Sigma^n(u)$ is

$$\langle, \rangle = -du^2 + \varrho^2(u)\langle, \rangle_\mathbb{P}. \tag{5.1}$$

It can be easily seen that a graph $\Sigma^n(u)$ is a spacelike hypersurface if and only if $|Du|_\mathbb{P}^2 < \varrho^2(u)$, $Du$ being the gradient of $u$ in $\mathbb{P}^n$ and $|Du|_\mathbb{P}$ its norm, both with respect to the metric $\langle, \rangle_\mathbb{P}$. The future-pointing Gauss map of a spacelike vertical graph $\Sigma^n(u)$ over $\Omega$ is given by the vector field

$$N(x) = \frac{\varrho(u(x))}{\sqrt{\varrho^2(u(x)) - |Du(x)|_\mathbb{P}^2}}\left(\partial_t|_{(u(x), x)} + \frac{1}{\varrho^2(u(x))}Du(x)\right), \quad x \in \Omega. \tag{5.2}$$

It is well known ([9, Lemma 3.1]) that in the case where $\mathbb{P}^n$ is a simply connected manifold, every complete spacelike hypersurface $\psi : \Sigma^n \to -I \times_\varrho \mathbb{P}^n$ such that the warping function $\varrho$ is bounded on $\Sigma^n$ is an entire spacelike graph over $\mathbb{P}^n$. In particular, this happens for complete spacelike hypersurfaces contained in a timelike bounded region of $-I \times_\varrho \mathbb{P}^n$. It is interesting to observe that, in contrast to the case of graphs into a Riemannian space, an entire spacelike graph $\Sigma^n(u)$ in a Lorentzian spacetime is not necessarily complete, in the sense that the induced Riemannian metric is not necessarily complete on $\mathbb{P}^n$. However, it can be proved that if $\mathbb{P}^n$ is complete and $|Du|_\mathbb{P}^2 \leq \varrho^2(u) - c$ for certain positive constant $c > 0$, then $\Sigma^n(u)$ is complete. Although a particular case of this claim is proven in [4, Theorem 4.1] and the general proof is analogous, we will expose it here for the sake of completeness.

**Proposition 1.** *Let $\mathbb{P}^n$ be a complete Riemannian manifold and $\Sigma^n(u)$ an entire spacelike vertical graph in $-I \times_\varrho \mathbb{P}^n$. If*
$$|Du|_\mathbb{P}^2 \leq \varrho^2(u) - c$$
*for certain positive constant $c > 0$, then $\Sigma^n(u)$ is complete.*

*Proof.* From (5.1), the Cauchy-Schwarz inequality and the assumptions of the proposition we get
$$\langle X, X \rangle = -\langle Du, X \rangle_\mathbb{P}^2 + \varrho^2(u) \langle X, X \rangle_\mathbb{P} \geq \left( \varrho^2(u) - |Du|_\mathbb{P}^2 \right) \langle X, X \rangle_\mathbb{P} \geq c \langle X, X \rangle_\mathbb{P},$$
for every $X \in \mathfrak{X}(\Sigma^n(u))$. This implies that $L \geq \sqrt{c} L_\mathbb{P}$, where $L$ and $L_\mathbb{P}$ denote the length of a curve on $\Sigma^n(u)$ with respect to the Riemannian metrics $\langle, \rangle$ and $\langle, \rangle_\mathbb{P}$, respectively. As a consequence, as $\mathbb{P}^n$ is complete by assumption, the induced metric on $\Sigma^n(u)$ from the metric of $-I \times_\varrho \mathbb{P}^n$ is also complete. □

In this context, we can establish a non-parametric version of Theorem 1.

**Theorem 4.** *Let $\mathbb{P}^n$ be a complete Riemannian manifold, and consider $-I \times_\varrho \mathbb{P}_f^n$ a spatially weighted GRW spacetime obeying (3.4) and such that the Hessian of the weight function $f$ is bounded from below. Let $\Sigma^n(u)$ be an entire vertical graph in $-I \times_\varrho \mathbb{P}_f^n$ determined by a bounded smooth function $u \in \mathcal{C}^\infty(\mathbb{P}^n)$. Suppose that the $f$-mean curvature $H_f$ of $\Sigma^n(u)$ satisfies*
$$(\log \varrho)'(u) \leq H_f \leq \alpha \quad \text{and} \quad H_f \geq 0,$$
*for some positive constant $\alpha$. If*

(5.3)
$$|Du|_\mathbb{P}^2 \leq \frac{\beta \inf_{\Sigma^n(u)} \varrho^2(u) \inf_{\Sigma^n(u)} |H_f - (\log \varrho)'(u)|^\gamma}{1 + \beta \inf_{\Sigma^n(u)} |H_f - (\log \varrho)'(u)|^\gamma}$$

*for some constants $\beta$ and $\gamma$, $\beta > 0$, $\gamma \neq 0$, then $u \equiv t_0$, for some $t_0 \in I$.*

*Proof.* Let us observe first that, under the assumptions of the theorem, $\Sigma^n(u)$ is complete. In fact, from (5.3) we easily obtain
$$\varrho^2(u) - |Du|_\mathbb{P}^2 \geq c = \frac{\inf_{\Sigma^n(u)} \varrho^2(u)}{1 + \beta \inf_{\Sigma^n(u)} |H_f - (\log \varrho)'(h)|^\gamma} > 0,$$
and the completeness of $\Sigma^n(u)$ follows from Proposition 1.

On the other hand, from (2.1), (2.3) and (5.2) we can get the following inequality
$$|\nabla h|^2 = \frac{|Du|_\mathbb{P}^2}{\varrho^2(u) - |Du|_\mathbb{P}^2}.$$
Therefore, the inequality (5.3) implies (4.2), and the result follows from Theorem 1. □

The non-parametric version of Theorem 2 can be stated and proved in a similar way as the previous one.

**Theorem 5.** *Let $\mathbb{P}^n$ be a complete Riemannian manifold, and consider $-I \times_\varrho \mathbb{P}_f^n$ a spatially weighted GRW spacetime obeying (3.4) with convex weight function $f$. Let $\Sigma^n(u)$ be an entire vertical graph in $-I \times_\varrho \mathbb{P}_f^n$ determined by a bounded smooth function $u \in \mathcal{C}^\infty(\mathbb{P}^n)$ and with constant $f$-mean curvature $H_f$ satisfying*
$$0 \leq H_f \sup_{\Sigma^n(u)} (\log \varrho)'(u) \leq H^2.$$
*If*
$$|Du|_\mathbb{P}^2 \leq \frac{\alpha \inf_{\Sigma^n(u)} \varrho^2(u) \left( \inf H^2 - H_f \sup_\Sigma (\log \varrho)'(u) \right)^\beta}{1 + \alpha \left( \inf H^2 - H_f \sup_\Sigma (\log \varrho)'(u) \right)^\beta}$$
*for some constants $\alpha$ and $\beta$, $\alpha > 0$, $\beta \neq 0$, then $u \equiv t_0$, for some $t_0 \in I$.*

In the case of Theorem 3 we need some extra assumptions in order to obtain a non-parametric version of it, since we cannot assure the completeness of an entire graph under the assumptions of Theorem 3.

**Theorem 6.** *Let $\mathbb{P}^n$ be a complete Riemannian manifold, and consider $-I \times_\varrho \mathbb{P}^n_f$ a spatially weighted GRW spacetime. Let $\Sigma^n(u)$ be an entire vertical graph in $-I \times_\varrho \mathbb{P}^n$ of a bounded smooth function $u \in C^\infty(\mathbb{P}^n)$ and suppose that $H_f > 0$, $\varrho'(u) > 0$ and that the following inequalities are satisfied*

$$\frac{n^2}{4}(\log \varrho)'^2(u) \leq \frac{n^2 H_f^2}{4} \leq (n-1)\frac{\varrho''}{\varrho}(u).$$

*If $|Du|^2_{\mathbb{P}} \leq \alpha \varrho^2(u)$, for some constant $0 < \alpha < 1$, and $\varrho(u) \in \mathcal{L}^p_f(\mathbb{P})$, for some $p > 1$, then $u \equiv t_0$, for some $t_0 \in I$, with $\mathrm{vol}_f(\Sigma^n(u)) < +\infty$. In addition, if $-I \times_\varrho \mathbb{P}^n$ obeys (3.4) and the weight function $f$ is bounded and convex, then $\Sigma^n(u)$ is compact.*

*Proof.* As we are assuming that $u$ is bounded and $|Du|^2_{\mathbb{P}} \leq \alpha \varrho^2(u)$ for some constant $0 < \alpha < 1$, it is easy to see that the assumptions of Proposition 1 are satisfied, so $\Sigma^n(u)$ is complete. Moreover, from equation (5.9) of [7] we get that

$$d\Sigma = \varrho^{n-1}(u)\sqrt{\varrho^2(u) - |Du|^2_{\mathbb{P}}}\, d\,\mathbb{P}.$$

Consequently, since $\varrho(u) \in \mathcal{L}^p_f(\mathbb{P})$, for $p > 1$, and $\varrho(h(q)) = \varrho(u(x))$ for all $q = (u(x), x) \in \Sigma^n(u)$, we obtain that $\varrho(h) \in \mathcal{L}^p_f(\Sigma^n(u))$, for $p > 1$. Therefore, we can apply Theorem 3 to conclude the result. □

## Acknowledgements

The first author is partially supported by Fundación Séneca project 19901/GERM/15, Spain and MINECO-FEDER Grant No MTM2015-65430-P. The second author is partially supported by CNPq, Brazil, grant 303977/2015-9. The third author is partially supported by CAPES, Brazil. The fourth author is partially supported by CNPq, Brazil, grant 308757/2015-7.

## References

1. R. Aiyama, *On the Gauss map of complete space-like hypersurfaces of constant mean curvature in Minkowski space*, Tsukuba J. Math. **16** (1992), 353–361.
2. A.L. Albujer and L.J. Alías, *Calabi-Bernstein results for maximal surfaces in Lorentzian product spaces*, J. Geom. Phys. **59** (2009), 620–631.
3. A.L. Albujer, F.E.C. Camargo and H.F. de Lima, *Complete spacelike hypersurfaces with constant mean curvature in $-\mathbb{R} \times \mathbb{H}^n$*, J. Math. Anal. Appl. **368** (2010), 650–657.
4. A.L. Albujer, F.E.C. Camargo and H.F. de Lima, *Complete spacelike hypersurfaces in a Robertson-Walker spacetime*, Math. Proc. Cambridge Philos. Soc. **151** (2011), 271–282.
5. J.A. Aledo and L.J. Alías, *On the curvatures of bounded complete spacelike hypersurfaces in the Lorentz-Minkowski space*, Manuscripta Math. **101** (2000), 401–413.
6. L.J. Alías and A.G. Colares, *Uniqueness of spacelike hypersurfaces with constant higher order mean curvature in generalized Robertson-Walker spacetimes*, Math. Proc. Cambridge Philos. Soc. **143** (2007), 703–729.
7. L.J. Alías, A.G. Colares and H.F. de Lima, *On the rigidity of complete spacelike hypersurfaces immersed in a generalized Robertson-Walker spacetime*, Bull. Braz. Math. Soc. **44** (2013), 195–217.
8. L.J. Alías, D. Impera and M. Rigoli, *Spacelike hypersurfaces of constant higher order mean curvature in generalized Robertson-Walker spacetimes*, Math. Proc. Cambridge Philos. Soc. **152** (2012), 365–383.
9. L.J. Alías, A. Romero and M. Sánchez, *Uniqueness of complete spacelike hypersurfaces of constant mean curvature in generalized Robertson-Walker spacetimes*, Gen. Relativity Gravitation **27** (1995), 71–84.
10. D. Bakry and M. Émery, *Diffusions hypercontractives*, In *Séminaire de probabilités*, XIX, 1983/84, volume 1123 of *Lecture Notes in Math.*, pages 177–206, Springer (Berlin, 1985).
11. G.P. Bessa, S. Pigola and A. Setti, *Spectral and stochastic properties of the $f$-Laplacian, solutions of PDEs at infinity and geometric applications*, Rev. Mat. Iberoam. **29** (2013), 579–610.
12. M. Caballero, A. Romero and R.M. Rubio, *Constant mean curvature spacelike surfaces in three-dimensional generalized Robertson-Walker spacetimes*, Lett. Math. Phys. **93** (2010), 85–105.
13. M. Caballero, A. Romero and R.M. Rubio, *Uniqueness of maximal surfaces in generalized Robertson-Walker spacetimes and Calabi-Bernstein type problems*, J. Geom. Phys. **60** (2010), 394–402.
14. E. Calabi, *Examples of Bernstein problems for some nonlinear equations*, Proc. Sympos. Pure Math. **15** (1970), 223–230.
15. F.E.C. Camargo, A. Caminha, H.F. de Lima and U. Parente, *Generalized maximum principles and the rigidity of complete spacelike hypersurfaces*, Math. Proc. Cambridge Philos. Soc. **153** (2012), 541–556.
16. J.S. Case, *Singularity theorems and the Lorentzian splitting theorem for the Bakry-Émery-Ricci tensor*, J. Geom. Phys. **60** (2010), 477–490.

17. M.P. Cavalcante, H.F. de Lima and M.S. Santos, *New Calabi-Bernstein type results in weighted generalized Robertson-Walker spacetimes*, Acta Math. Hungar. **145** (2015), 440–454.
18. S.Y. Cheng and S.T. Yau, *Maximal space-like hypersurfaces in the Lorentz-Minkowski spaces*, Ann. of Math. **104** (1976), 407–419.
19. H.F. de Lima, A.M.S. Oliveira and M.S. Santos, *Rigidity of complete spacelike hypersurfaces with constant weighted mean curvature*, Beitr. Algebra Geom., Online First. DOI: 10.1007/s13366-015-0253-7.
20. H.F. de Lima and M.A.L. Velásquez, *Uniqueness of complete spacelike hypersurfaces via their higher order mean curvatures in a conformally stationary spacetime*, Math. Nachr. **287** (2014), 1223–1240.
21. A. Grigor'yan and L. Saloff-Coste, *Dirichlet heat-kernel in the exterior of a compact set*, Comm. Pure Appl. Math. **55** (2002), 93–133.
22. M. Gromov, *Isoperimetry of waists and concentration of maps*, Geom. Funct. Anal. **13** (2003), 178–215.
23. J.M. Latorre and A. Romero, *Uniqueness of noncompact spacelike hypersurfaces of constant mean curvature in generalized Robertson-Walker spacetimes*, Geom. Dedicata **93** (2002), 1–10.
24. G. Li and I.M.C. Salavessa, *Graphic Bernstein results in curved pseudo-Riemannian manifolds*, J. Geom. Phys. **59** (2009), 1306–1313
25. B. O'Neill, *Semi-Riemannian Geometry with Applications to Relativity*, Academic Press (New York, 1983).
26. S. Pigola, M. Rigoli and A.G. Setti, *Vanishing theorems on Riemannian manifolds, and geometric applications*, J. Funct. Anal. **229** (2005), 424–461.
27. M. Rimoldi, *Rigidity results for Lichnerowicz Bakry-Émery Ricci tensors*, Ph.D. thesis, Università degli Studi di Milano (Milano, 2011).
28. A. Romero and R.M. Rubio, *On the mean curvature of spacelike surfaces in certain three-dimensional Robertson-Walker spacetimes and Calabi-Bernstein's type problems*, Ann. Global Anal. Geom. **37** (2010), 21–31.
29. A. Romero, R. Rubio and J.J. Salamanca, *Uniqueness of complete maximal hypersurfaces in spatially parabolic generalized Robertson-Walker spacetimes*, Classical Quantum Gravity **30** (2013), 115007 pp. 13.
30. A. Romero, R. Rubio and J.J. Salamanca, *A new approach for uniqueness of complete maximal hypersurfaces in spatially parabolic GRW spacetimes*, J. Math. Anal. Appl. **419** (2014), 355–372.
31. G. Wei and W. Wylie, *Comparison geometry for the Bakry-Émery-Ricci tensor*, J. Differential Geom. **83** (2009), 377–405.
32. Y.L. Xin, *On the Gauss image of a spacelike hypersurface with constant mean curvature in Minkowski space*, Comment. Math. Helv. **66** (1991), 590–598.

[1] Departamento de Matemáticas, Campus Universitario de Rabanales, Universidad de Córdoba, 14071 Córdoba, Spain.

*E-mail address*: `alma.albujer@uco.es`

[2] Departamento de Matemática, Universidade Federal de Campina Grande, 58429-970 Campina Grande, Paraíba, Brazil.

*E-mail address*: `henrique@mat.ufcg.edu.br`
*E-mail address*: `arlandsonm@gmail.com`
*E-mail address*: `marco.velasquez@pq.cnpq.br`